\documentclass[11pt]{article}
\usepackage{amsmath,amssymb}
\usepackage{epsfig}
\begin{document}

 \vspace*{1\baselineskip} \begin{center}\Large\bf
Complexity in linearly coupled dynamical networks: Some unusual
phenomena in energy accumulation \footnote {\small This work is
supported by the National
 Science Foundation of
 China under grants 60674093, 60334030, 10472001. }
 \end{center}
\vspace*{1\baselineskip}

\centerline {Zhisheng Duan, ~Jinzhi Wang, ~Guanrong Chen ~and ~Lin
Huang}

\vspace*{0.5\baselineskip}
\begin{center}
 The  State Key Laboratory for Turbulence and Complex Systems and
 Department of Mechanics and
       Engineering Science, Peking University,  Beijing  100871, P. R.
       China. \\  Emails:
       duanzs@pku.edu.cn,  eegchen@cityu.edu.hk, hl35hj75@pku.edu.cn  \\
       {\it Tel and Fax}: (8610)62765037
\end{center}
\vskip 0.5cm

{\bf Abstract.} \,\, This paper addresses the energy accumulation
problem, in terms of the $H_2$ norm, of linearly coupled dynamical
networks. An interesting outer-coupling relationship is constructed,
under which the $H_2$ norm of the newly constructed network with
column-input and row-output shaped matrices increases exponentially
fast with the node number $N$: it increases generally much faster
than $2^N$ when $N$ is large while the $H_2$ norm of each node is 1.
However, the $H_2$ norm of the network with a diffusive coupling is
equal to $\gamma_2 N$, i.e., increasing linearly, when the network
is stable, where $\gamma_2$ is the $H_2$ norm of a single node. And
the $H_2$ norm of the network with antisymmetrical coupling also
increases, but rather slowly, with the node number $N$. Other
networks with block-diagonal-input and block-diagonal-output
matrices behave
 similarly. It demonstrates that the changes of $H_2$ norms
in different networks are very complicated, despite the fact that
the networks are linear. Finally, the influence of the $H_2$ norm of
the locally linearized network on the output of a network with Lur'e
nodes is discussed.

{\bf Keywords.} \,\,$H_2$ norm of linear network, Diffusive
coupling, Antisymmetrical coupling.

\section{ Introduction}

Complex  networks have attracted increasing attention from
physicists, biologists, social scientists and control engineers in
recent years [1-4, 9, 11-13]. Many complicated problems in various
networks, such as topological structures, small-world and scale-free
characteristics, robustness and fragility, and self-similarity, have
been studied (see [1, 4, 11-13] and references therein). In general,
 a complex dynamical network can be considered as a large-scale
system with special interconnections among its dynamical nodes. The
large-scale system theory has been extensively studied in the last
three decades, and many interesting results have been established,
on such basic issues as decentrally fixed modes, decentralized
controllers design, diagonal Lyapunov function method, and M-matrix
method \cite{sil91}. Many complicated technical problems in
large-scale systems, such as decomposition and aggregation,
connective stability and instability, competitive equilibria, local
instability, have also been studied (see \cite{duan05, duan06,
sil91} and references therein). New concepts such as `small-world',
'scale-free', `power law', etc. in the current studies of complex
networks will in turn provide new opportunities to the traditional
large-scale system theory. On the other hand, the mature large-scale
system methods can  provide effective tools to network stability
analysis and control problems.
 Note that many dynamical aspects of networks such as network synchronization,
network control, and dynamical modeling have  been studied (see
\cite{bar04, bel05, li04, wang02}). From a control theoretic
viewpoint, a dynamical system  has inputs and outputs, so does a
network, where the input and output relation has to be considered.
 This paper attempts to explore
the input to output $H_2$ gain (i.e., $H_2$ norm, which represents
the energy of the system output with fixed input) of a general
linear dynamical network with different couplings. The results show
that the $H_2$-norm changes can be unexpectively complicated,
despite the linear structure of the network.

The rest of this paper is organized as follows. In section 2, some
preliminaries are introduced. In section 3, an interesting network
with column-input and row-output matrices is constructed, whose
energy increases exponentially fast with the node number $N$.
   In section 4,  diffusively  and
antisymmetrically coupled networks are studied. Compared with the
network studied in section 3, the $H_2$ norms of these  network
increases linearly and in particular the norm of the diffusive
network is equal to the node number $N$. In section 5, networks with
block-diagonal-input and block-diagonal-output matrices are
similarly discussed. In section 6, for a Lur'e network, the
influence of the energy of its locally linearized network on the
original network output is studied.  The last section concludes the
paper.

\section{Preliminaries}

Consider a  continuous-time linear system,
\begin{equation} \label{n1}
\dot{\xi}=A_1\xi+B_1w_1,\quad y_1=C_1\xi, \end{equation} where $\xi$
is the state of
 the system, $w_1$ and $y_1$ are the input and output of the system, respectively, and
  $A_1\in {\bf R}^{n\times
n}, B_1\in {\bf R}^{n\times m}, C_1\in {\bf R}^{l\times n}$ are
given constant matrices. The transfer function from $w_1$ to $y_1$
is $G_1(s)=C_1(sI-A_1)^{-1}B_1$. If $A_1$ is stable, the $H_2$ norm
of system (\ref{n1}) is represented by the $H_2$ norm of the
transfer function $G_1(s)$, which is defined by
$$\|G_1(s)\|_2=\sqrt{\frac{1}{2\pi}\int_{-\infty}^{+\infty}\hbox{trace}\{G_1^H(jw)G_1(jw)\}dw}.$$
The $H_2$ norm represents the power of the system output for fixed
input, or the induced system gain, in control theory \cite{zhou96}.
For computing the $H_2$ norm, the following formula is convenient
 \cite{zhou96}.

{\bf Lemma 1}\,\, If $A_1$ is stable, then the $H_2$ norm of system
(\ref{n1}) is given by
$$\|G_1(s)\|_2=\sqrt{\hbox{trace}(B_1^TQB_1)},$$
where $Q$ satisfies the Lyapunov equation
$$QA_1+A_1^TQ+C_1^TC_1=0.$$  \hfill $\Box$

When $A_1$ is unstable, the $L_2$ norm can be computed
\cite{zhou96}. With a linear coupling, the $N$ nodes, each described
by (\ref{n1}), constitute a dynamical network as follows:
 \begin{equation}  \label{sy1}
\dot{x}_i=A_1x_i+\sum_{\tiny
\begin{array}{c}j=1\end{array}}^N \gamma_{ij} A_{12} x_j,
\quad i=1, \cdots, N,
\end{equation}
 where $\gamma_{ij} \in {\bf R}, $ $A_1$ is given as in (\ref{n1}), and $A_{12}\in {\bf R}^{n\times
n}$ is the inner coupling matrix describing the interconnections
among components of $x_j, j=1,\dots,N$. Let
$\Gamma=(\gamma_{ij})_{N\times N}$, where
 $\gamma_{ij}$ are given in (\ref{sy1}), which is referred to as the
outer coupling matrix. Using the Kronecker product notation
\cite{fax04, li04}, a network with column-input and row-output
shaped matrices can be rewritten as
 \begin{equation}  \label{sy2}
\dot{x}=(I_N\bigotimes A_1+\Gamma\bigotimes A_{12})x+B w_1, \quad
y=C x, \end{equation} where $x=(x_1^T, \cdots, x_N^T)^T$,
$B=E_N\bigotimes B_1$, $C=E_N^T\bigotimes C_1$  and $E_N=(1, \cdots,
1)^T\in {\bf R}^N$.

To discuss the $H_2$ norm of network (\ref{sy2}), the following
simple result for stability analysis is needed.

{\bf Lemma 2}\,\, Suppose $\lambda_i, i=1,\cdots, k$, are
 distinct eigenvalues of $\Gamma$, with $k\leq N$, not to count the multiplicity. Then,  $I_N\bigotimes
A_1+\Gamma\bigotimes A_{12}$ is stable if, and only if,
$A_1+\lambda_iA_{12}, i=1, \cdots, k$, are stable simultaneously.

{\bf Proof}\,\, Let $T$ be a nonsingular matrix such that
$T^{-1}\Gamma T=J$, where $J$ is the Jordan form of $\Gamma$. Then,
the following similarity transformation completes the proof easily:
$$ (T^{-1}\bigotimes I)(I_N\bigotimes A_1+\Gamma\bigotimes A_{12})(T\bigotimes
I)=I_N\bigotimes A_1+J \bigotimes A_{12}.$$ \hfill $\Box$

Throughout this paper, let $A_{12}=B_1C_1$, which simply means that
the network has an input and output inner coupling.

\section{ $H_2$-norm energy accumulation }

In this section, consider the changes of $H_2$ norm with the network
size for a specially constructed network. For convenience, first
consider the following $N\times N$ coupling matrix:
\begin{equation} \label{ga}
\Gamma=(\gamma_{ij})_{N\times N}, \, \gamma_{ii}=0; \,
\gamma_{ij}=-1, \, \hbox{if}\,\, i>j, \, j\leq i-2 \,\,
\hbox{and}\,\, j \,\, \hbox{is odd; otherwise},\,\, \gamma_{ij}=1.
\end{equation}
One can also define $\Gamma$ inductively, that is,
\begin{equation} \label{ga2}
\Gamma=\left(\begin{array}{cc} \Gamma_{(N-1)\times (N-1)} & \gamma_1\\
\gamma_2 & 0\end{array}\right),\,\, N\geq 2, \end{equation}
 where $\gamma_1$ is a column vector with  all elements being 1, $\gamma_2$
 is a row vector with  the $i$-th element being  -1 if $i\leq N-2$  is odd and the other
 elements being 1.  For
example, when $N=5$, one has
$$\Gamma_{5\times 5}=\left(\begin{array}{ccccc} 0 & 1 & 1 & 1 & 1\\1
& 0 & 1 & 1 & 1\\ -1 & 1 & 0 & 1 & 1\\-1 & 1 & 1 & 0 & 1\\ -1 & 1 &
-1 & 1 & 0\end{array}\right).$$

{\bf Lemma 3}\,\, When $N$ is even, the distinct eigenvalues of
$\Gamma$ are -1 and 1, and $\Gamma$ is similar
to $$\Lambda_N=\left(\begin{array}{ccc} -I_{\frac{N}{2}} & 0 \\
0 & J_{\frac{N}{2}} \end{array}\right),$$ where $I_{\frac{N}{2}}$ is
the identity matrix of order $\frac{N}{2}\times \frac{N}{2}$ and
$J_{\frac{N}{2}}$ is a Jordan matrix of order $\frac{N}{2}\times
\frac{N}{2}$ as
 $$J_{\frac{N}{2}}=\left(\begin{array}{cccc} 1 & 2 & &\\
 & \ddots &\ddots & \\& & \ddots &  2\\ & & &
 1\end{array}\right).$$
 When $N$ is odd, the distinct eigenvalues of
$\Gamma_{N\times N}$ are -1, 1 and 0, and  $\Gamma_{N\times N}$ is
similar to $\hbox{diag} (\Lambda_{N-1}, 0)$.

{\bf Proof} \,\,It suffices to prove this lemma by  inductively
determining the eigenvectors and generalized eigenvectors of
$\Gamma$.

First, suppose that $N$ is even. If $N=2$, obviously, the
eigenvalues of $\Gamma_{2\times 2}$ are -1 and 1. If $N=4$, one has
$$ \Gamma_{4\times 4}P_2=-P_2 \quad \hbox{and} \quad \Gamma_{4\times
4}Q_2=Q_2J_2, $$ where $P_2=\left(\begin{array}{cccc}0 & -1 & 0 &  1\\
 0 & 0 & 1 & -1 \end{array}\right)^T$,
$Q_2=\left(\begin{array}{cccc} 1& 1& 0 & 0\\-1 & -1 &1 &
1\end{array}\right)^T $ and  $J_2$ is a Jordan matrix of order 2 as
defined above. This shows that the statement holds for $N=4$. Let
$P_{\frac{N-2}{2}}$ and $P_{\frac{N}{2}}$ ($N\geq 6$) be
$(N-2)\times \frac{N-2}{2} $ and $N \times \frac{N}{2} $ matrices
composing of column eigenvectors of $\Gamma_{(N-2)\times (N-2)}$ and
$\Gamma_{N\times N}$ corresponding to the eigenvalue -1,
respectively. Then, $P_{\frac{N}{2}}$ is composed of $(0,\,-1,\,
0,\,1,\, 0_{1\times(N-4)})^T$ and columns $(0_{\frac{N-2}{2}\times
2}, \, P^T_{\frac{N-2}{2}})^T$.

Similarly, for $N\geq 6$,
 with $Q_2$ as given above, let $Q_{\frac{N-2}{2}}$ and $Q_{\frac{N}{2}}$  be $
(N-2)\times\frac{N-2}{2}$ and $N\times \frac{N}{2} $ matrices
composing of column and generalized column eigenvectors of
$\Gamma_{(N-2)\times (N-2)}$ and $\Gamma_{N\times N}$ corresponding
to the eigenvalue 1, respectively. Then, $Q_{\frac{N}{2}}$ is
composed of $\left(\begin{array}{c}Q_{\frac{N-2}{2}} \\
0_{2\times \frac{N-2}{2}}\end{array}\right)$ and
$q_{\frac{N}{2}}=(-2\alpha, -2\alpha, \beta_1, \cdots, \beta_{N-4},
1, 1)^T$, where $\alpha$ is the first element of the last column of
$Q_{\frac{N-2}{2}}$ and $\beta_i$ in $q_{\frac{N}{2}}$, $i=1,\cdots,
N-4$, are determined by
 $$(\Gamma_{N\times N}-I_{N\times N})q_{\frac{N}{2}}=2(q^T_{\frac{N-2}{2}}, 0, 0)^T, $$
where  $q_{\frac{N-2}{2}}$ is the last column of
$Q_{\frac{N-2}{2}}$.

 When $N$ is odd, except a zero eigenvalue, it is similar to the
 case $N$ being even. \hfill $\Box$

{\bf Remark 1}\,\,Column eigenvectors for the eigenvalue -1 and
generalized column  eigenvectors for the eigenvalue 1 are discussed
in the proof of Lemma 3. Row and generalized row eigenvectors for
the eigenvalues -1 and 1 can be studied similarly.

 With the coupled matrix as discussed above,
the stability of network (\ref{sy2}) can be easily tested.

 {\bf Corollary 1}\,\, If the number of nodes, $N$,
is even, then network (\ref{sy2}) is stable, i.e.,  $A=I \bigotimes
A_1 + \Gamma \bigotimes A_{12}$ is stable if, and only if,
$A_1-A_{12}$ and $A_1+A_{12}$ are stable. If $N$ is odd, $A=I
\bigotimes A_1 + \Gamma \bigotimes A_{12}$ is stable if, and only
if, $A_1,$ $A_1-A_{12}$ and $A_1+A_{12}$ are stable. \hfill $\Box$

{\bf Theorem 1}\,\, If the number of nodes, $N$, is even, the $H_2$
norm of system (\ref{sy2}) is equal to the $H_2$ norm of the
following system:
 \begin{equation}   \label{th21}
 \dot{\xi}=(I\bigotimes A_1+ J_{\frac{N}{2}}\bigotimes A_{12})\xi+B_{\frac{N}{2}}w,\quad y=C_{\frac{N}{2}}\xi,
\end{equation}
where $J_{\frac{N}{2}}$ is given as in Lemma 1, and
$C_{\frac{N}{2}}$ and $B_{\frac{N}{2}}$ are determined inductively
as follows:
$$C_{\frac{N}{2}}=(C_{\frac{N}{2}-1}, \,
(-1)^{\frac{N}{2}-1}2^{\frac{N}{2}-2}C_1), \,\, \hbox{with}\,\,
C_2=(2C_1,\,0),$$ in which, by denoting
$B_{\frac{N}{2}-1}=(2B_1^T,\,x_1B_1^T,\,x_2B_1^T,\,\cdots,\,x_{\frac{N}{2}-3}B_1^T,\,
B_1^T)^T$,
$$B_{\frac{N}{2}}=(2B_1^T,\,(2+x_1)B_1^T,\,(x_1+x_2)B_1^T,\,\cdots,\,(x_{\frac{N}{2}-3}+1)B_1^T,\,
B_1^T)^T,\,\, \hbox{with}\,\, B_2=(2B_1^T,\,B_1^T)^T. $$ If $N(N\geq
5)$ is odd, the $H_2$ norm of system (\ref{sy2}) is equal to the
$H_2$ norm of the following system:
 \begin{equation}   \label{th22}
 \dot{\xi}=\tilde{A}_o\xi+\tilde{B}_{\frac{N+1}{2}o}w,\quad y=\tilde{C}_o\xi,
\end{equation}
where $\tilde{A}_o=\hbox{diag}(I\bigotimes A_1+
J_{\frac{N-1}{2}}\bigotimes A_{12}, \, 0)$,
$\tilde{C}_o=(C_{\frac{N-1}{2}}, \, (-1)^{\frac{N-1}{2}}C_1),$ and
 $\tilde{B}_{\frac{N+1}{2}o}$ is determined inductively
as follows:
$B_{\frac{N-1}{2}o}=(2B_1^T,\,x_1B_1^T,\,x_2B_1^T,\,\cdots,\,x_{\frac{N-1}{2}-2}B_1^T,\,
2B_1^T,\,B_1^T)^T$,  and {\small
$$B_{\frac{N+1}{2}o}=(2B_1^T,\,(2+x_1)B_1^T,\,(x_1+x_2)B_1^T,\,\cdots,\,(x_{\frac{N+1}{2}-2}+2)B_1^T,\,
2B_1^T,\,B_1^T)^T,\,\, \hbox{with}\,\,
B_{2o}=(2B_1^T,\,2B_1^T,\,B_1^T)^T. $$}

{\bf Proof} \,\, Taking a similarity transformation to change
$\Gamma$ into a Jordan form as in Lemma 3, one can complete the
proof easily. \hfill $\Box$

{\bf Remark 2}\,\, From the Jordan matrix $J_{\frac{N}{2}}$, one can
see that the number 2 appears in it. This means that the term $2^N$
will appear in the transfer function of system (\ref{th21}) or
(\ref{th22}). It is imaginable that  the $H_2$ norm of network
(\ref{sy2}) can be possibly very large. Of course, one may also take
other similarity transformations such that the term 2 becomes 1 in
$J_{\frac{N}{2}}$. However, it should be noted that  any similarity
transformation does not change the system norm.

{\bf Example 1}\,\,
 Consider the single node system (\ref{n1}) with
      $A_1=-4.5,\, B_1=1,\, C_1=3.$
See Table 1 for the $H_2$  norms of network (\ref{sy2}).

 {\bf Example
2}\,\, Consider the single node system (\ref{n1}) with
$$ A_1=\left(\begin{array}{cc} 0 & 1\\-4 &
-2.5\end{array}\right),\quad
 B_1=\left(\begin{array}{c} 0 \\ 1\end{array}\right),\quad
 C_1=\left(\begin{array}{cc}  2& 2\end{array}\right).$$
See Table 1 for the $H_2$ norms of network (\ref{sy2}).

\begin{center}
{\small Table 1 \,\, $H_2$ norms of network  (\ref{sy2}) with
different $N$.} \vskip 3mm $\begin{tabular}{|c|c|c|c|c|c|c|c|c|c|c|}
\hline $ N$ & $1$ & $ 2$ &
 $3$& $4$ &$5$ & $\cdots $ & $10 $ & $\cdots $& $20 $& $\cdots$\\
\hline $\hbox{Ex 1}$ & $1 $ &
$3.4641 $ & $6.0828 $&$15.4919$& $25.4755$ & $\cdots $ & $1.4363\times 10^3 $ & $\cdots $& $3.6864\times 10^6 $& $\cdots$\\
\hline
 $\hbox{Ex 2}$ & $1 $ &
$4.8990 $ & $9.1761 $& $40.9878 $ &$74.9765$ & $\cdots $&
$3.6676\times 10^4 $& $\cdots$ & $4.2018\times 10^9$
& $\cdots$ \\
\hline
\end{tabular}$
\end{center}

From Table 1, one can see that the $H_2$ norm of network (\ref{sy2})
increases exponentially fast, much faster than $2^{20}$ ($=
1048576$) when $N$ is large. This means that even a single node
system can generate a complicated network through simple coupling
relations.

\section{Diffusive and antisymmetrical coupling relationships}

First, consider  diffusive coupling networks. Change the outer
coupling matrix $\Gamma$ in (\ref{ga}) as follows:
 \begin{equation} \label{gad}
  \Gamma_D=(\gamma_{ij})_{N\times N}, \,\hbox{if}\, i\not=j, \, \gamma_{ij} \,
  \hbox{are given as in (\ref{ga})} \,\hbox{and}\, \gamma_{ii}=-
  \sum_{j=1, j\not=i}^N \gamma_{ij}, \, i,j=1,\cdots,N.
  \end{equation}
Obviously, the sum of every row of $\Gamma_D$ is zero, which
generally refers to as diffusive coupling.

{\bf Lemma 4}\,\, If $N$ is odd, the eigenvalues of $\Gamma_D$,
given as in (\ref{gad}),  are $-N, -(N-2), -(N-2), \cdots, -3, -3,
-1, 0$; especially, if $N=3$, its eigenvalues are -3, -1 and 0. If
$N$ is even, the eigenvalues of $\Gamma_D$ are $-N, -(N-2), -(N-2),
\cdots, -2, -2, 0$; especially, if $N=2$, its eigenvalues are -2 and
0. \hfill $\Box$

{\bf Corollary 2}\,\, Given $A_1$ and $A_{12}$ with compatible
dimensions, $I_N\bigotimes A_1+\Gamma_D\bigotimes A_{12}$ is stable
if, and only if, $A_1, \, A_1-A_{12},\, A_1-3A_{12},\, \cdots, \,
A_1-(N-2)A_{12}$ and $A_1-NA_{12}$ are all stable when $N$ is odd,
or $A_1, \, A_1-2A_{12},\, \cdots, \, A_1-(N-2)A_{12}$ and
$A_1-NA_{12}$ are all stable when $N$ is even. \hfill $\Box$

{\bf Theorem 2}\,\, If network (\ref{sy2}) with $\Gamma=\Gamma_D$ is
stable, then the $H_2$ norm of (\ref{sy2}) is $\gamma_2 N$, where
$\gamma_2$ is the $H_2$ norm of  the single-node system (\ref{n1}).

{\bf Proof}\,\,Because of the diffusive characteristic of
$\Gamma_D$, 0 is one of its  eigenvalues and the corresponding
eigenvector is ${\bf 1}_N=(1,\, \cdots, \, 1)^T$. Let  $P$ be a
nonsingular matrix such that $P^{-1}\Gamma_DP=J$, where $J$ is the
Jordan form of $\Gamma_D$, and the last column of $P$ is ${\bf 1}_N$
and the last row of $J$ is a zero vector. Let
$\eta=(P^{-1}\bigotimes I)x$. Then, systems (\ref{sy2}) becomes
 \begin{equation} \label{th31}
 \dot{\eta}=(I_N\bigotimes A_1+J\bigotimes A_{12}) \eta+ (P^{-1}{\bf
 1}_N \bigotimes B_1)w,\,\, y=({\bf 1}_N^TP\bigotimes C_1) \eta.
 \end{equation}
 Obviously, $P^{-1}{\bf  1}_N=(0,\,\cdots,\,0,\,\alpha)^T,\,
 \alpha\not=0$. Let the last element of ${\bf 1}_N^TP$ be $\beta$.
 Then, one has
 $$\alpha\beta = {\bf 1}_N^TP P^{-1}{\bf  1}_N =N.$$
 Therefore, the transfer function from $w$ to $y$ is $N
 C_1(sI-A_1)^{-1}B_1$. This completes the proof. \hfill $\Box$

{\bf Corollary 3}\,\, For any outer coupling matrix $\Gamma$, if the
sum of every its row is 0, or the sum of every its column is 0, then
the $H_2$ norm of  network (\ref{sy2}) is $\gamma_2 N$ when it is
stable, where $\gamma_2$ is the $H_2$ norm of the single-node system
(\ref{n1}). \hfill $\Box$

{\bf Remark 3}\,\, By Theorem 2 or Corollary 3, the $H_2$ norm of
system (\ref{sy2}) with $\Gamma_D$ given in (\ref{gad}) or with any
diffusive matrix is equal to $20$ when $N=20$ in Examples 1 and 2
discussed above, verifying the linear growth speed of the network
energy.

Next, consider networks with an antisymmetrical coupling matrix,
 \begin{equation} \label{gaa}
  \Gamma_A=(\gamma_{ij})_{N\times N}, \,\hbox{if}\, i>j, \, \gamma_{ij}=-1; \,
  \hbox{if} \, i<j, \, \gamma_{ij}=1; \, \hbox{and} \, \gamma_{ii}=0, \, i,j=1,\cdots,N.
  \end{equation}
 Obviously, all the eigenvalues of $\Gamma_A$ are located on the imaginary
 axis. Table 2 shows the $H_2$ norms of network (\ref{sy2}) with
 the outer coupling matrix $\Gamma=\Gamma_A$ given by (\ref{gaa}).

\begin{center}
{\small Table 2 \,\, $H_2$ norms of network  (\ref{sy2}) with
$\Gamma=\Gamma_A$ for system data given in Section 3.} \vskip 3mm
$\begin{tabular}{|c|c|c|c|c|c|c|c|c|c|c|} \hline $ N$ & $1$ & $ 2$ &
 $3$& $4$ &$5$ & $\cdots $ & $10 $ & $\cdots $& $20 $& $\cdots$\\
\hline $\hbox{Ex 1}$ & $1 $ &
$1.8397 $ & $2.5542 $&$3.1978$& $3.8013$ & $\cdots $ & $6.6168 $ & $\cdots $& $12.1951 $& $\cdots$\\
\hline
 $\hbox{Ex 2}$ & $1 $ &
$1.7949 $ & $2.4851 $& $3.1381 $ &$3.7805$ & $\cdots $& $7.0102$&
$\cdots$ & $13.5857$
& $\cdots$ \\
\hline
\end{tabular}$
\end{center}

Similarly to the diffusive coupling case, one can see from Table 2
that the $H_2$ norm of network (\ref{sy2}) with  an antisymmetrical
coupling also increases slowly.

\section{Networks with  block-diagonal-input and block-diagonal-output matrices}

 In this section, consider network (\ref{sy2}) with block-diagonal-input and
 block-diagonal-output matrices, i.e.,
  \begin{equation} \label{m1}
   B=\hbox{diag}(B_1,\cdots,B_1),\, C=\hbox{diag}(C_1,\cdots,C_1),
 \end{equation}
but $A$ remains unchanged. Table 3 shows the $H_2$ norm changes of
network (\ref{sy2}) with (\ref{m1}).

\begin{center}
{\small Table 3 \,\, $H_2$ norms of network  (\ref{sy2}) with
(\ref{m1}) and  $\Gamma$ given  in (\ref{ga}) for system data given
in Section 3.} \vskip 3mm $\begin{tabular}{|c|c|c|c|c|c|c|c|c|c|c|}
\hline $ N$ & $1$ & $ 2$ &
 $3$& $4$ &$5$ & $\cdots $ & $10 $ & $\cdots $& $20 $& $\cdots$\\
\hline $\hbox{Ex 1}$ & $1 $ &
$1.8974 $ & $2.7203 $&$6.9089$& $11.4222$ & $\cdots $ & $1.0026\times 10^3 $ & $\cdots $& $3.7236\times 10^6 $& $\cdots$\\
\hline
 $\hbox{Ex 2}$ & $1 $ &
$2.5531 $ & $4.3147 $& $23.5579 $ &$43.9595$ & $\cdots $&
$3.4090\times 10^4$& $\cdots$ & $5.5723\times 10^9$
& $\cdots$ \\
\hline
\end{tabular}$
\end{center}

Comparing Table 3 with Table 1, one can see that the norms in Table
3 also increase dramatically. For small $N$, the norms in Table 3
are smaller than the ones in Table 1. However, they become larger
than the ones in Table 1 when $N$ is large.

\begin{center}
{\small Table 4 \,\, $H_2$ norms of network (\ref{sy2}) with
(\ref{m1}) and $\Gamma_D$ given  in (\ref{gad}) for the same system
data.} \vskip 3mm $\begin{tabular}{|c|c|c|c|c|c|c|c|c|c|c|} \hline $
N$ & $1$ & $ 2$ &
 $3$& $4$ &$5$ & $\cdots $ & $10 $ & $\cdots $& $20 $& $\cdots$\\
\hline $\hbox{Ex 1}$ & $1 $ &
$1.1952 $ & $1.4918 $&$1.5780$& $1.7967$ & $\cdots $ & $2.2347 $ & $\cdots $& $2.9020 $& $\cdots$\\
\hline
 $\hbox{Ex 2}$ & $1 $ &
$1.1602 $ & $1.4915 $& $1.5494 $ &$1.8108$ & $\cdots $& $2.2498$&
$\cdots$ & $2.9944$
& $\cdots$ \\
\hline
\end{tabular}$
\end{center}

\begin{center}
{\small Table 5 \,\, $H_2$ norms of network (\ref{sy2}) with
(\ref{m1}) and $\Gamma_A$ given as in (\ref{gaa}) for the same
system data.} \vskip 3mm $\begin{tabular}{|c|c|c|c|c|c|c|c|c|c|c|}
\hline $ N$ & $1$ & $ 2$ &
 $3$& $4$ &$5$ & $\cdots $ & $10 $ & $\cdots $& $20 $& $\cdots$\\
\hline $\hbox{Ex 1}$ & $1 $ &
$1.4142 $ & $1.7321 $&$2$& $2.2361$ & $\cdots $ & $3.1623 $ & $\cdots $& $4.4721 $& $\cdots$\\
\hline
 $\hbox{Ex 2}$ & $1 $ &
$1.4591 $ & $1.8108 $& $2.1012 $ &$2.3533$ & $\cdots $& $3.3311$&
$\cdots$ & $4.7108$
& $\cdots$ \\
\hline
\end{tabular}$
\end{center}

Comparing Tables 4 and 5 with Theorem 2 and Table 2, one can see
that the $H_2$ norms of network (\ref{sy2}) with (\ref{m1}) also
increase slowly for the couplings $\Gamma_D$ and $\Gamma_A$, even
much slower than the norms of the networks with column-input and
row-output matrices.

{\bf Corollary 4} For first-order node systems, the $H_2$ norm of
its corresponding network with antisymmetrical coupling matrix
$\Gamma_A$ is $\gamma_2 \sqrt{N}$.

{\bf Proof}\,\, By Lemma 1, in this case, $Q$ can be taken as a
diagonal matrix in the form of $\alpha I$, and $B^TQB=\gamma_2^2N$,
where $\gamma_2$ is the $H_2$ norm of the single-node system. This
completes the proof. \hfill $\Box$

{\bf Remark 4}\,\, One can clearly verify the result of Corollary 4
by examining Table 5.

\section{Influence of energy on nonlinear Lur'e networks}

The absolute stability of Lur'e systems and synchronization of Lur'e
networks have been extensively studied (see  \cite{cur97, wu02} and
references therein). In this section, consider the influence of the
energy of the locally linerized system of a nonlinear Lur'e network
on its network output. Given a Lur'e system as follows:
\begin{equation} \label{L1}
\left\{ \begin{array}{l} \dot{x}_1=(A_1-2A_{12})x_1+B_{01}f_1(y_1), \quad x_1(0)=B_{01}, \\
 y_1=C_{01}x_1, \end{array} \right. \end{equation}
 where $x_1$ is the state, $y_1$ is the measured output,
 $x_1(0)$ is the initial condition,  $$ A_1=\left(\begin{array}{cc} 0 & 1\\-4 &
-2.5\end{array}\right),\quad
 B_{01}=\left(\begin{array}{c} 0 \\ 1\end{array}\right),\quad
C_{01}=\left(\begin{array}{cc}  2& 2\end{array}\right),\quad
A_{12}=B_{01}C_{01},$$ referring to Example 2 in Section 3 for these
matrix data, and the nonlinear function $f_1(y_1)=|y_1+1|-|y_1-1|$,
which is just a piece wise linear function as used in the canonical
Chua's circuit \cite{wu02}. Obviously, $f_1$ satisfies  the
following sector condition:
 \begin{equation} \label{sec12}
0\leq \frac{f_{1}(y_{1})}{y_{1}}\leq 2, \quad
f_{1}(0)=0.\end{equation} Consider a network with regulated output
formed by the node equation (\ref{L1}), as follows:
\begin{equation} \label{NL1}
\left\{\begin{array}{l} \dot{x}=(I_N\bigotimes (A_1-2A_{12})+\Gamma_L\bigotimes A_{12})x+Bf(y),\quad x(0)=E_N \bigotimes x_1(0), \\
 y=C_1x,\\z=C_2x, \end{array} \right.  \end{equation}
 where $x=(x_1^T, \cdots, x_N^T)^T$, $x_i$ is the state of the $i$-th Lur'e
 system, $y=(y_1, \cdots, y_N)^T$,  $f(y)=(f_1(y_1), \cdots, f_1(y_N))^T$,
$B=I_N \bigotimes B_{01},$ $C_1=I_N \bigotimes C_{01},$ $C_2=E^T_N
\bigotimes C_{01},$ $z$ is
 the regulated output,
   $E_N$ is given as in (\ref{sy2}), and $\Gamma_L$ is the outer
   coupling matrix similarly to $\Gamma$ in (\ref{sy2}).

  Linearizing network (\ref{NL1}) near the zero equilibrium, one gets a
  linear network as follows:
 \begin{equation} \label{NL2}
\left\{\begin{array}{l} \dot{x}=(I_N\bigotimes A_1+\Gamma_L\bigotimes A_{12})x,\quad x(0)=E_N \bigotimes x_1(0), \\
 z=C_2x. \end{array} \right.  \end{equation}
By the LQR control method \cite{zhou96}, network (\ref{NL2}) can be
viewed in another form as follows:
 \begin{equation} \label{NL3}
\left\{\begin{array}{l} \dot{x}=(I_N\bigotimes A_1+\Gamma_L\bigotimes A_{12})x+E_N \bigotimes x_1(0)\delta(t),\quad x(0_-)=0, \\
 z=C_2x, \end{array} \right.  \end{equation}
where $\delta(t)$ is the impulse function.
 It is well-known that the output energy of network (\ref{NL2}) is
 equal to the $H_2$ norm of network (\ref{NL3}), i.e.,
 $$\|C_2(sI-I_N\bigotimes A_1-\Gamma_L\bigotimes A_{12})^{-1}E_N \bigotimes
 x_1(0)\|_2.$$
According to the discussions in Sections 2 and 3, with the matrix
data in (\ref{L1}), the output energy of (\ref{NL2}) increases
exponentially fast with $\Gamma_L=\Gamma$ as in (\ref{ga}) and
increases linearly with $\Gamma_L=\Gamma_D$ as in (\ref{gad}). Fig.
1 and Fig. 2 show the influence of energy changes of network
(\ref{NL2}) or (\ref{NL3}) on the regulated output of Lur'e network
(\ref{NL1}).
 \noindent
\begin{center}
 \unitlength=1cm
 \hbox{\hspace*{0cm} \epsfxsize5cm \epsfysize5cm
\epsffile{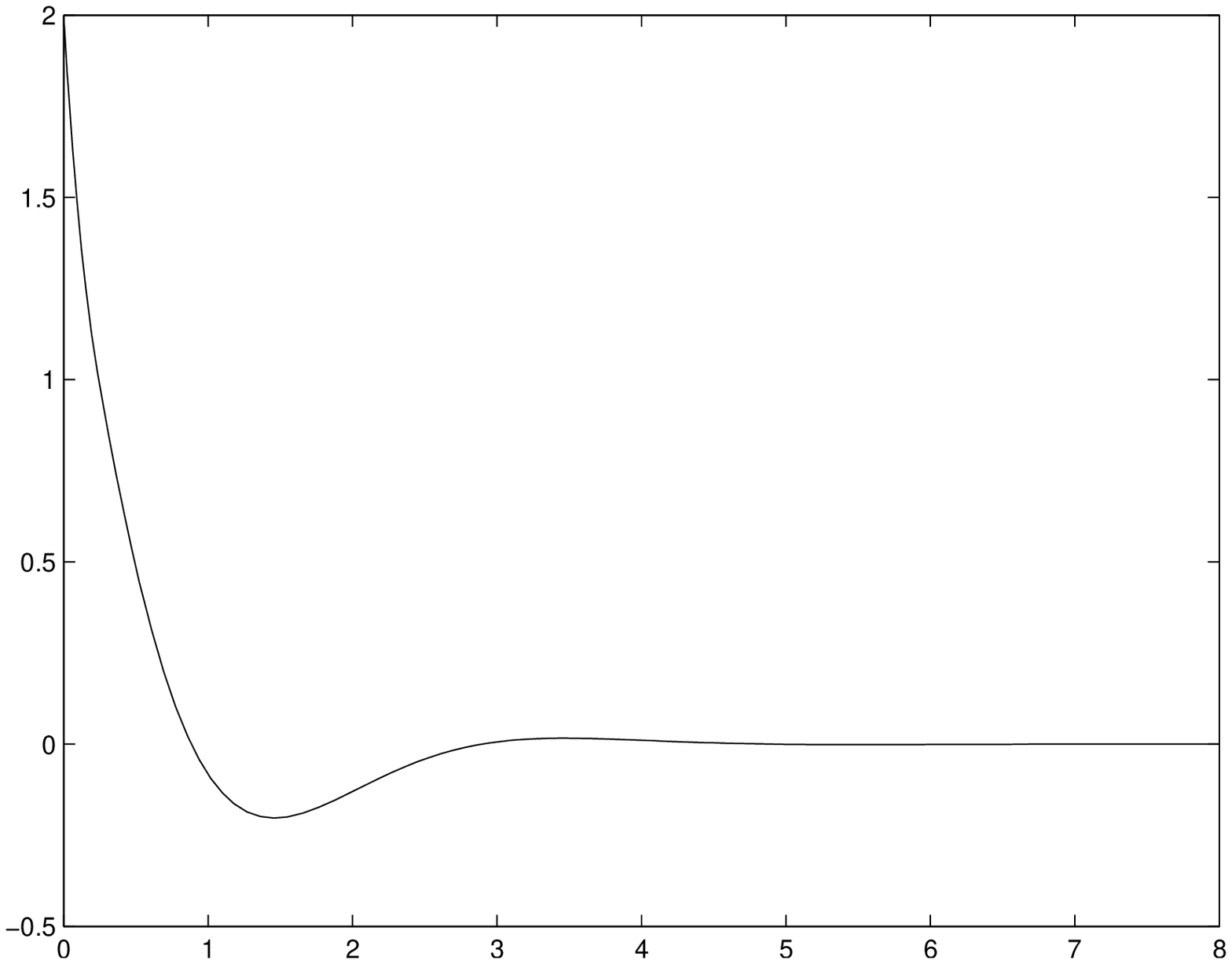}\quad \epsfxsize5cm \epsfysize5cm
\epsffile{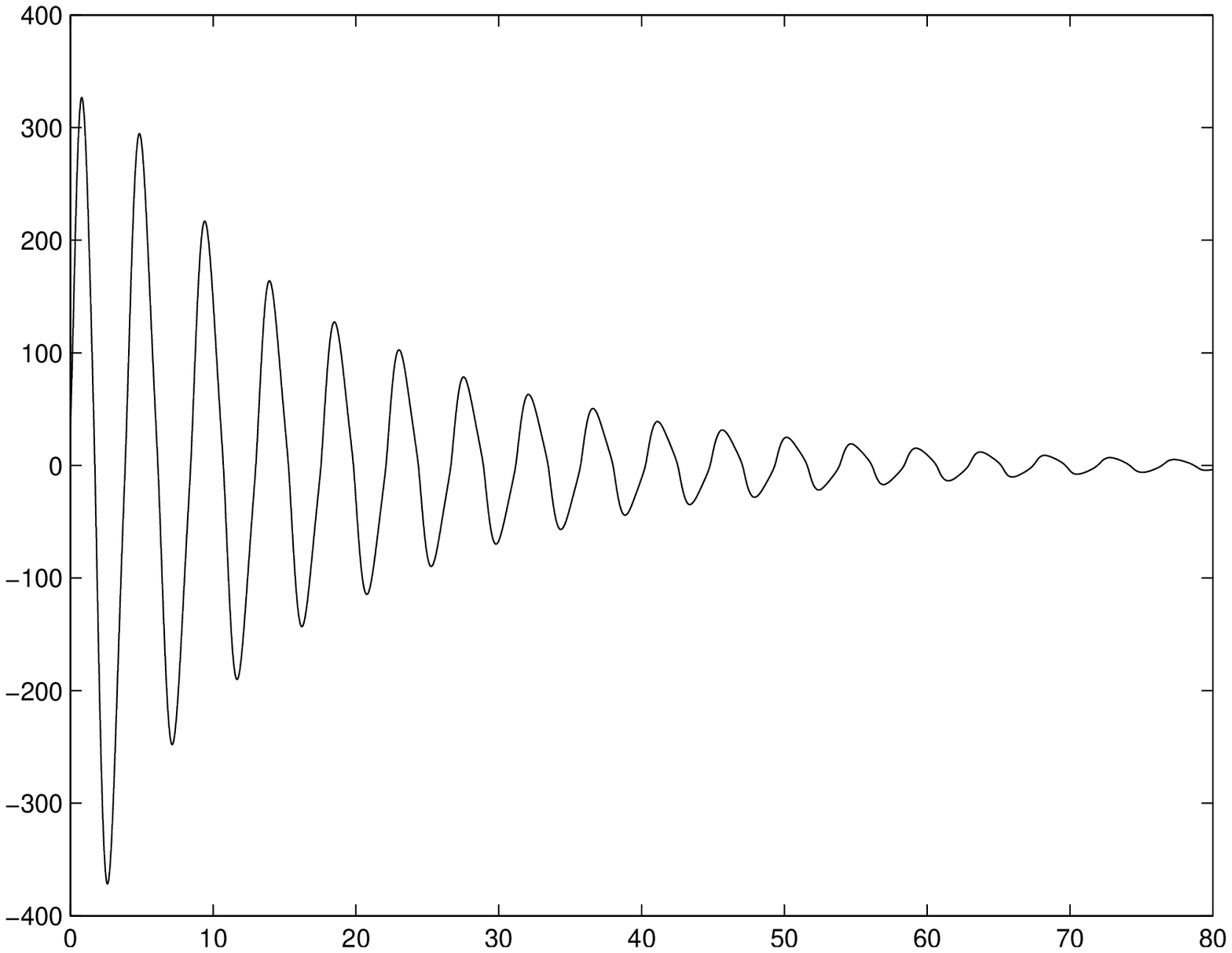}\quad \epsfxsize5cm \epsfysize5cm
\epsffile{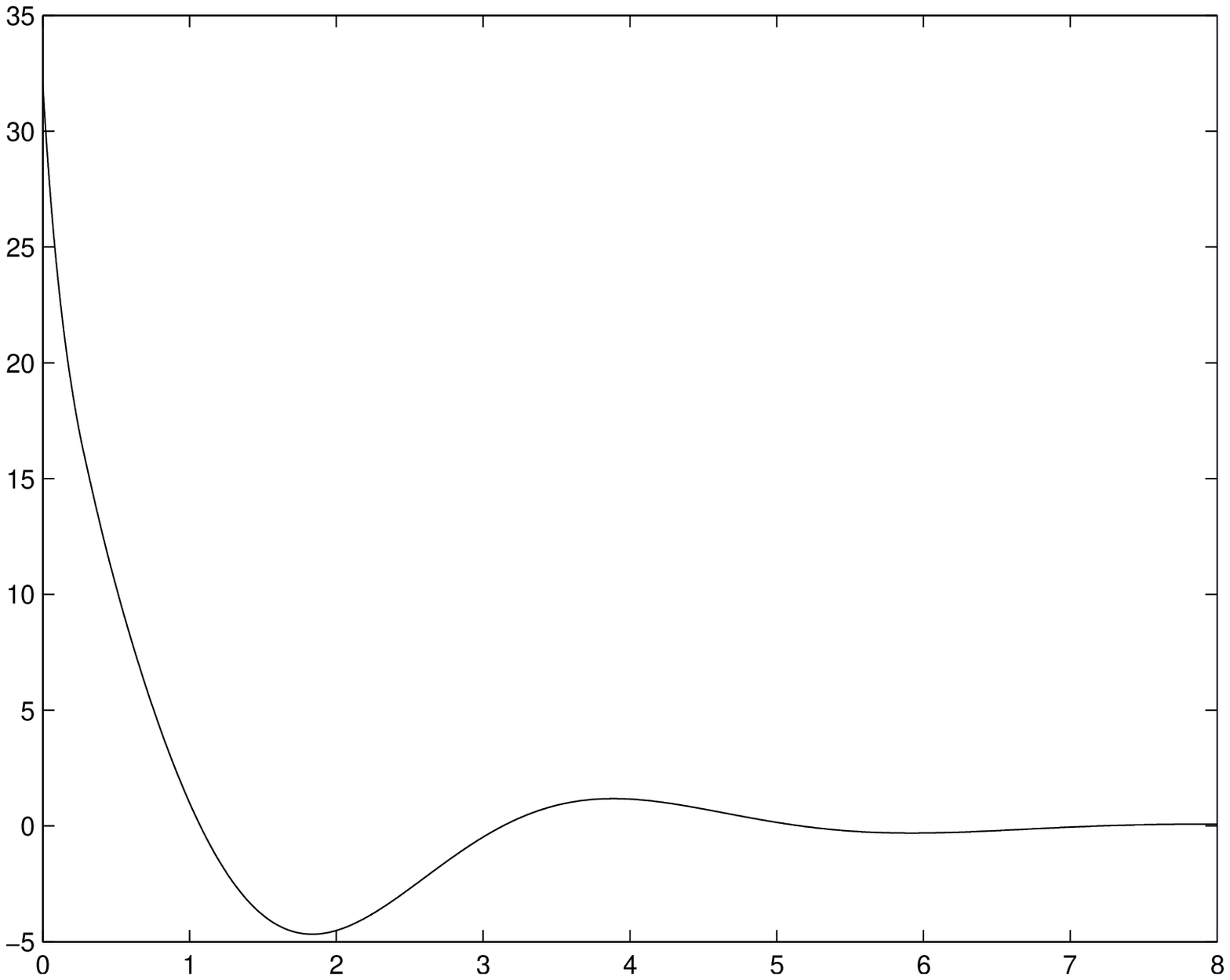}
 }
\end{center}
\vskip -1.0cm
 \centerline{\quad(a) The output of Lur'e \qquad (b) The output of network (\ref{NL1}) \quad (c) The output of network (\ref{NL1})}
 \centerline{\quad  system (\ref{L1}). \qquad\qquad\qquad\quad with $\Gamma_L=\Gamma$ as in (\ref{ga}).
 \qquad\qquad \qquad with $\Gamma_L=\Gamma_D$ as in (\ref{gad}).}

 \centerline{ Fig. 1 \,\, The output of a single-node system (\ref{L1}) and
  outputs of Lur'e networks for $N=16$.}

\begin{center}
 \unitlength=1cm
 \qquad \hbox{\hspace*{0.1cm}  \epsfxsize5.5cm \epsfysize5cm
\epsffile{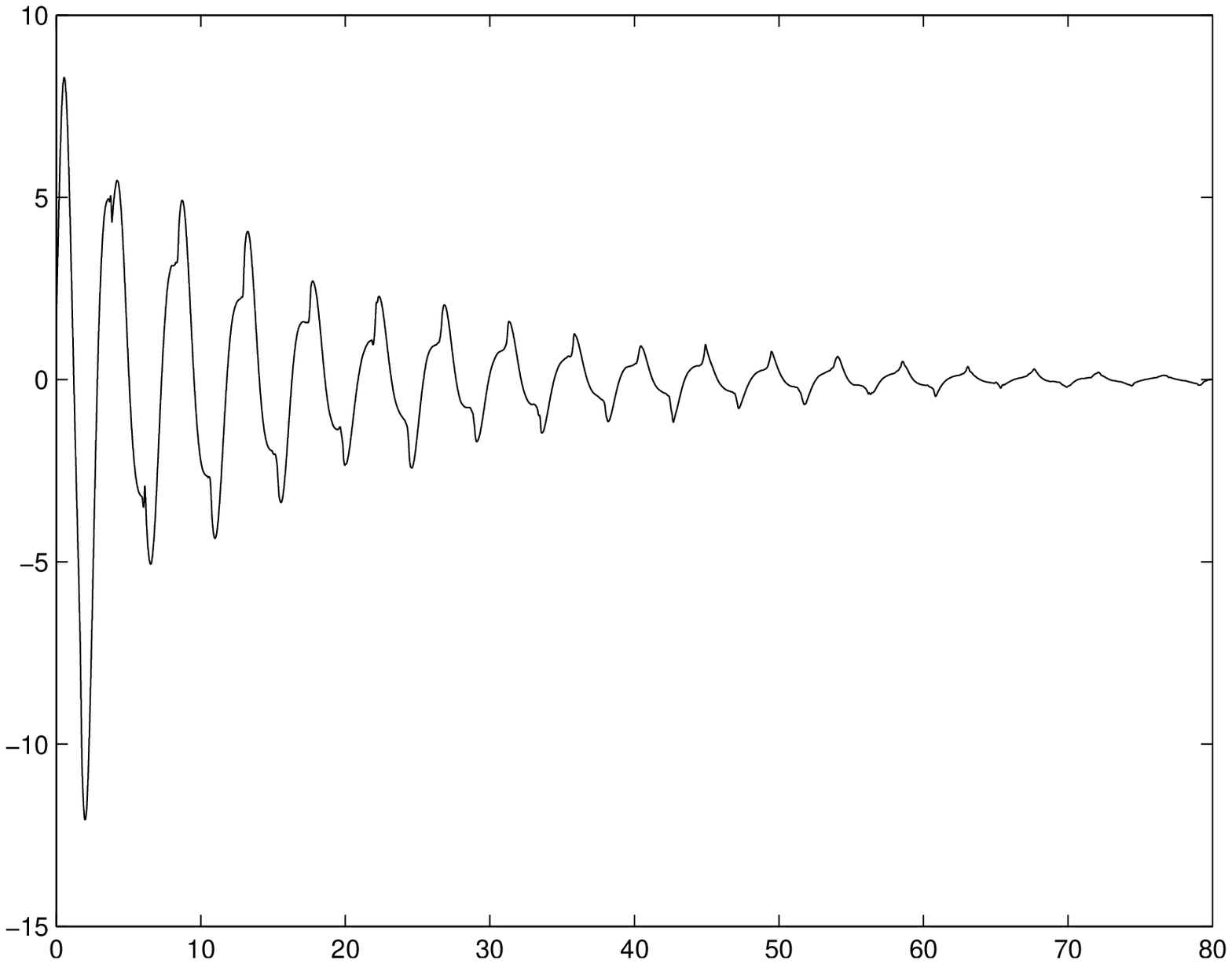}\qquad \epsfxsize5.5cm \epsfysize5cm
\epsffile{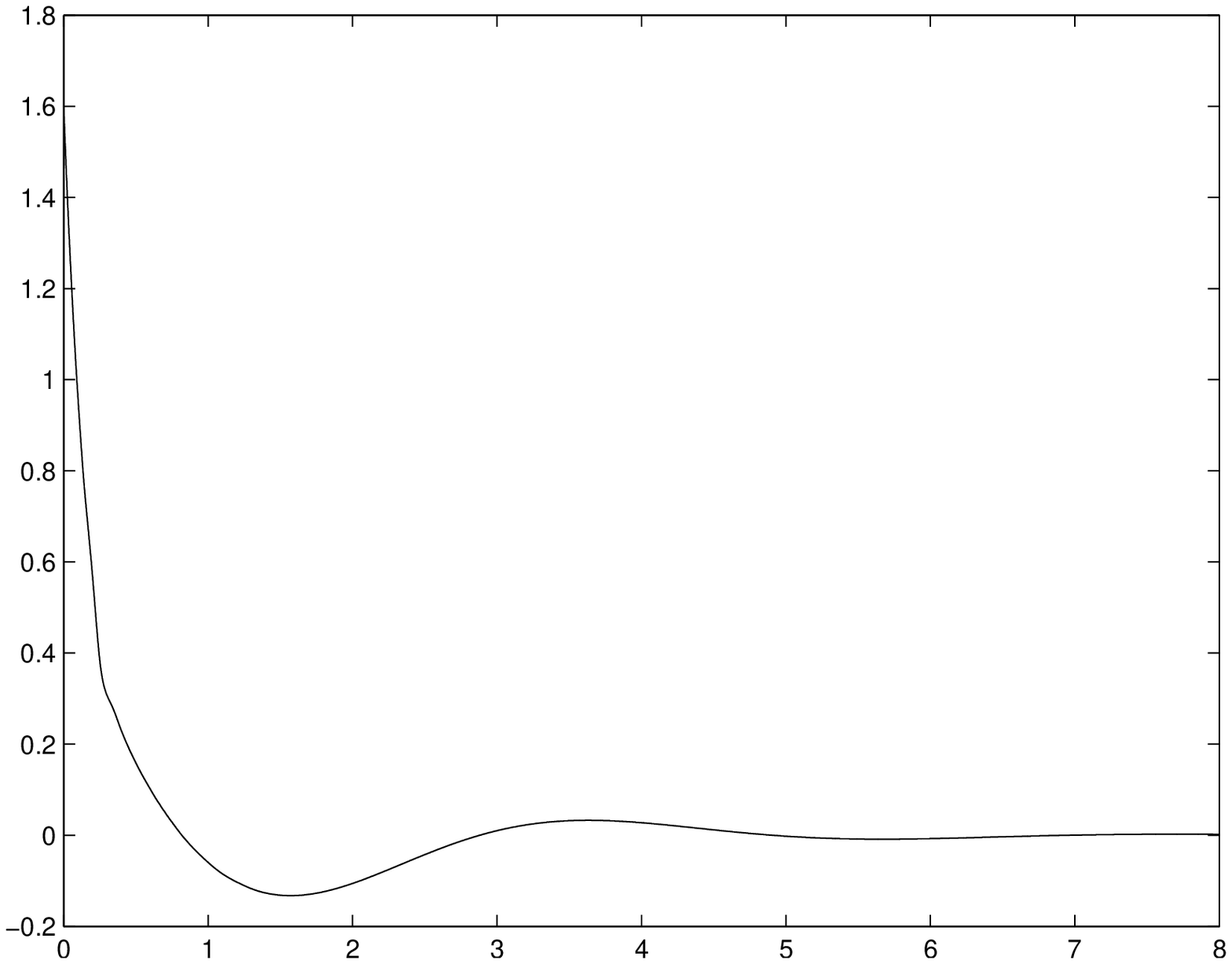}
 }
\end{center}
\vskip -0.3cm \centerline{\qquad (a) The output error of network
(\ref{NL1}) \quad (b) The output error of network (\ref{NL1})}
 \centerline{\qquad with $\Gamma_L=\Gamma$ as in (\ref{ga}).
 \qquad  \qquad\qquad with $\Gamma_L=\Gamma_D$ as in (\ref{gad}).}
 \centerline{ Fig. 2 \,\, The output errors of Lur'e networks with two initial
 conditions:}
 \centerline{ $x(0)=E_N \bigotimes x_1(0)$ and  $x(0)=0.95E_N \bigotimes x_1(0)$, for $N=16$.}

From Fig. 1 and Fig. 2, one can see that the shape of the solution
of the diffusive network is  similar to the solution of a
single-node Lur'e system, but the output error of the diffusive
network is smaller, for the two initial conditions shown in Fig. 2
(b). However, for the newly constructed coupling relationship
$\Gamma_L=\Gamma$ as in (\ref{ga}), the network output is very
different with the output of a single-node Lur'e system and the
network output error with two initial conditions is larger.

{\bf Remark 5}\,\, In fact, by the above idea, an input to output
performance index can be established for network synchronization
problems, especially for local synchronization problems. It is
imaginable that the synchronization performance is bad for a
nonlinear network when the energy of its locally linearized network
is large. Therefore,  large energy changes should be avoided for
synchronization problems.

\section{Conclusion}

In this paper, the $H_2$-norm energy accumulation problem has been
addressed for linearly coupled dynamical networks. Three types of
networks, i.e., a newly constructed network, a typical diffusively
coupled network and an antisymmetrically coupled network, have been
studied on the changes of their $H_2$ norms with respect to  the
network size $N$. For the newly constructed network, the $H_2$ norm
increases exponentially fast, even much faster than $2^N$ when $N$
is large. However, the $H_2$ norms of the diffusively coupled and
antisymmetrically coupled  networks increase only linearly or even
slower. This shows the complexity in energy changes of dynamical
networks, despite the linear dynamical nature of the networks, which
may be unexpected or even surprising. It has also been shown that
large energy changes can have bad influence on network output,
therefore should be avoided in network synchronization problems.
Finally, it should be noted that the corresponding $H_\infty$ norm
problem \cite{chui97} can be similarly studied, which will be
reported elsewhere.

 \hspace*{34pt}

\end{document}